\documentclass[11pt]{article}
\usepackage[utf8]{inputenc}
\usepackage{filecontents}
\usepackage[many]{tcolorbox}
\usepackage{graphicx}
\usepackage{harpoon}
\usepackage{xcolor}
\usepackage{latexsym}
\usepackage{amssymb}
\usepackage{amsmath}
\usepackage{mdframed}
\usepackage{hyperref, url}
\usepackage[show]{ed}
\usepackage{subfigure}
\usepackage[ntheorem]{empheq} 
\usepackage{amsthm, amssymb, amsfonts, latexsym}
\usepackage{pst-all}
\mathchardef\mhyphen="2D

\makeatletter
\DeclareRobustCommand{\cev}[1]{%
  \mathpalette\do@cev{#1}%
}
\newcommand{\do@cev}[2]{%
  \fix@cev{#1}{+}%
  \reflectbox{$\m@th#1\vec{\reflectbox{$\fix@cev{#1}{-}\m@th#1#2\fix@cev{#1}{+}$}}$}%
  \fix@cev{#1}{-}%
}
\newcommand{\fix@cev}[2]{%
  \ifx#1\displaystyle
    \mkern#20mu
  \else
    \ifx#1\textstyle
      \mkern#20mu
    \else
      \ifx#1\scriptstyle
        \mkern#26mu
      \else
        \mkern#26mu
      \fi
    \fi
  \fi
}

\makeatother

\newcommand{\fig}[1]{\epsfbox}
\newcommand{\bp}{\begin{minipage}{3.1cm}}
\newcommand{\ep}{\end{minipage}}

\newtheorem{Theorem}{Theorem}[]

\makeatletter 
\renewcommand{\thefigure}{\@arabic\c@figure}
\makeatother
\usepackage{xr}

\RequirePackage{graphicx,xcolor}
\RequirePackage{colortbl}
\RequirePackage{booktabs}
\RequirePackage{algorithm}
\RequirePackage[noend]{algpseudocode}
\RequirePackage{changepage}
\RequirePackage[twoside,%
				letterpaper,includeheadfoot,%
				layoutsize={8.125in,10.875in},%
                layouthoffset=0.1875in,%
                layoutvoffset=0.0625in,%
                left=38.5pt,%
                right=43pt,%
                top=43pt,% 10pt provided by headsep
                bottom=32pt,%
                headheight=0pt,% No Header
                headsep=10pt,%
                footskip=25pt,
                marginparwidth=38pt]{geometry}
\RequirePackage[labelfont={bf,sf},%
                labelsep=period,%
                figurename=Fig. ]{caption}
\DeclareCaptionFormat{pnasformat}{\normalfont\sffamily\fontsize{7}{9}\selectfont#1#2#3}
\definecolor{lightgray}{gray}{0.95}

\newtcolorbox{paperbox}[1][]{
    enhanced,
    colback=white,
    boxrule=0pt,
    boxsep=0pt,
    arc=0mm,
    width=0.8\linewidth,
    fuzzy shadow={0mm}{-4pt}{-4pt}{1mm}{black!30!white},
    #1
}

\begin{document}

\title{Universal set of Observables for Forecasting Physical Systems through Causal Embedding}

\author{G Manjunath*, A de Clercq* \& MJ Steynberg$^\dag$\\
{\footnotesize \it *Department of Mathematics \& Applied Mathematics, $^\dag$~Department of Physics, University of Pretoria, Pretoria 0028} \\ {\small Email: manjunath.gandhi@up.ac.za ; declercq@uchicago.edu ; thys.steynberg@fau.de} \\}
\date{}
\date{}
\maketitle

%\begin{multicols}{2}

\begin{abstract}
 We show how a pair of points can uniquely represent a left-infinite sequence obtained from observations of an underlying dynamical system through a phenomenon called  \textit{causal embedding}. A driven dynamical system creates such pairs, and a function can be learnt on them that can reconstruct the underlying dynamics as in Takens delay embedding.  The approach assures embedding stability unlike Takens delay embedding, and learnability which can be absent while stability can be present in the current reservoir computing framework. This accurately models underlying systems where recent methods like the next-generation reservoir computing fail. We demonstrate results and compare with the other methods including SINDY-PI.
 \end{abstract}


\section{Introduction}
Performing experimental measurements on  physical, biological and artificial deterministic systems to obtain a more informative dynamical model with applications in better control and management of resources and service interruptions has been well established in modern-day science. Attempts to build models from complex data such as sunspot cycles date back to the 1920s in \cite{yule1927vii}. Although our ability to collect data has vastly improved ever since, and data-driven machine learning methods (e.g., \cite{jaeger2004harnessing, schmid2010dynamic,williams2015data, xu2006modeling,brunton2016discovering,lu2018attractor, champion2019data, kaheman2020sindy, gauthier2021next,grigoryeva2023learning, fuhg2024extreme}) have gained popularity, the quest for the right observables, like the delay-coordinates in Takens delay embedding \cite{takens1981detecting} along with stability in forecasting has remained elusive.

In a customary forecasting problem of a dynamical system, one has an underlying or an unknown dynamical system $(W,S)$ where $S : W \to W$ is an invertible map, and one has access only to observations $(\theta(w_i), \theta(w_{i+1}), \ldots)$ of an orbit $(w_i, w_{i+1}, \ldots)$, where $w_{n+1}=S(w_n)$ and $\theta: W \to \mathbb{R}$ is an observation. In general there is no function from $\theta(w_i) \mapsto w_{i+1}$. 

The Takens delay embedding  theorem and its generalizations  (e.g., \cite{sauer1991embedology,stark1999delay,gutman2018embedding}) essentially guarantee that when observations when $W$ is a space whose dimension can be well-defined
(like a manifold), and has dimension $d$, the functions $S$, and $\theta: W \to \mathbb{R}$  belong to a generic class of smooth functions, then the delay coordinates  $\Phi_{2d,\theta}:W \to \mathbb{R}^{2d+1}$  given by 
            $\Phi_{2d,\theta}(w) := (\theta(S^{-2d}w)\ldots,\theta(S^{-1}w),\theta(w))$  is an embedding of $W$ in $\mathbb{R}^{2d+1}$ (i.e., the map $\Phi_{2d,\theta}: W \to \Phi_{2d,\theta}(W)$ is a homeomorphism onto its image).  Consequently, there exists a 
            function/map between the delay coordinates
 $F_S: (\theta(S^{-2d}(w)),\ldots,\theta(S^{-1}(w)),\theta(w)) \mapsto 
            (\theta(S^{-2d+1}(w)),\ldots,\theta(w),\theta(S(w))$. Hence, by observing successive data points in the space of delay coordinates, one can learn the map $F_S$ which then leads to forecasting of the observations obtained from $\theta$.  The salient feature of Takens delay embedding is that the delay coordinate map by $\Phi_{2d,\theta}$ is a universal observable of the data, i.e., for a fixed delay $d$, it does not change (but only adjusts its domain) even if the underlying $W$, $S$ or $\theta$ are changed. 
Next, Takens' result makes $\Phi_{2d,\theta}$ a conjugacy (we refer to \cite{de2013elements} for basic terminologies) that results in $F_S \circ \Phi_{2d,\theta} = \Phi_{2d,\theta} \circ S$ (or the diagram in Fig~\ref{Fig_COMM}A in green commutes).
Since the conjugacy establishes a one-to-one correspondence between the orbits of $S$ and $F_S$, a model $F_S$ is particularly helpful in understanding the global behavior of $S$ when
            long-term pointwise prediction is not possible, for instance, when $S$ exhibits sensitive dependence on initial conditions.

\begin{figure}	
 \centering
\includegraphics[width=9cm]{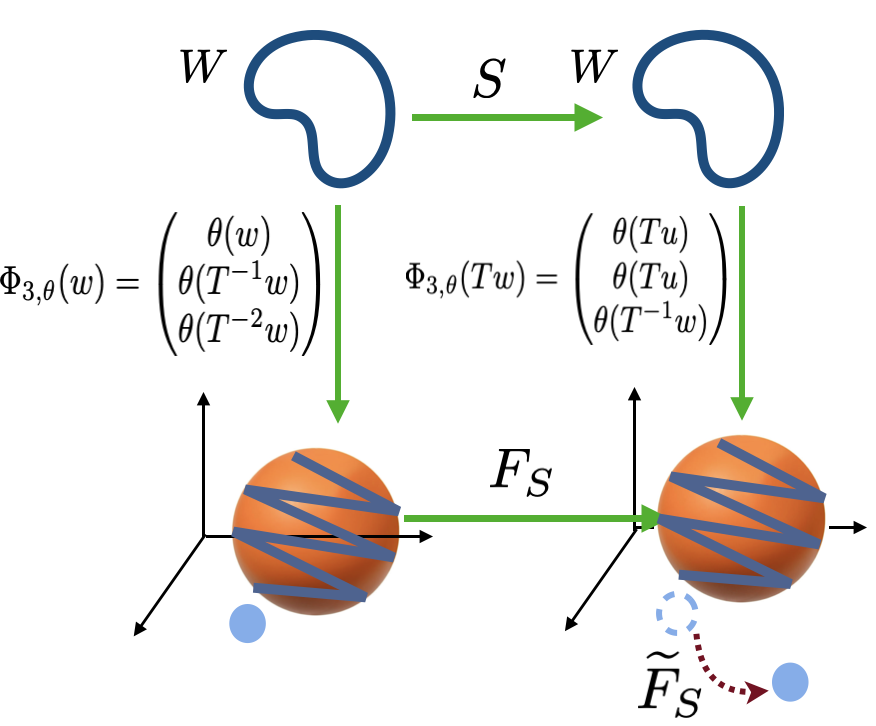}
\caption{Schematic of employing the Takens embedding, and the effect of iterates of $\widetilde{F}_S$  slipping out of $\Phi_{\theta,2d}(W)$.} 
\label{Fig_Big_A}
\end{figure}

The above results all hold when a delay-coordinate map $\Phi_{K,\theta}$ is used when $K \ge 2d+1$. 
In practice, iterating $F_S$ in the presence of errors can amplify the error depending on $\theta$ and $K$ used (see \cite{casdagli1991state} for a detailed discussion). In particular, we note
   the delay-coordinates can take values outside the embedded set $\Phi_{2d,\theta}(W)$ (see Fig.~\ref{Fig_Big_A}A), and hence 
         $\widetilde{F}_S$, an extension \cite{munkres2014topology} of the map $F_S$ that is defined on a neighborhood of $\Phi_{2d,\theta}(W)$ 
         is used to generate the dynamics in practice (e.g., by a feedforward neural network). While doing so,  the $F_S$-embedded object $\Phi_{2d,\theta}(W)$ 
          cannot be expected to be an attractor of an extension $\widetilde{F}_S$ since the extension  $\widetilde{F}_S$ behaves much like $F_S$ in a neighborhood of $\Phi_{2d,\theta}(W)$. Thus $\widetilde{F}_S$ 
         can be locally expanding in a neighborhood of $\Phi_{2d,\theta}(W)$. Hence for an $x$ lying outside $\Phi_{2d,\theta}(W)$, its future iterates  could be veering far off from $\Phi_{2d,\theta}(W)$ thereby leading to miscarried forecasting/ modeling (see Fig~\ref{Fig_Big_A}).

The problem of trajectories veering off is a stability issue and the resultant of the constructed dynamics on the set of interest not showing dissipative dynamics. More precisely, we may not have an attractor of $\widetilde{F}_S$ that contains $\Phi_{2d,\theta}(W)$. The problem has not been overcome even with the recent advances in learning dynamical systems. Recent Koopman operator-based data-driven methods use non-standardized observable selection (e.g., \cite{schmid2010dynamic,williams2015data}).  Library-based methods \cite{xu2006modeling,brunton2016discovering,champion2019data, kaheman2020sindy}) that require the vector field or a map to lie in the span of a pre-determined set of basis functions give interpretable models, however there is no guarantee that the vector field generated has an attractor. We demonstrate the unstable phenomenon on a more sophisticated version  of vector field algorithm in SINDY-PI  \cite{kaheman2020sindy}.  We consider the complete coordinates generated from a double pendulum (see Equations \eqref{DP_Eqn1m} and \eqref{DP_Eqn2m} in the Appendix), and show that while in the short-term, the generated system of differential equations generate a trajectory close-by to the original one, but in the long run, the trajectories diverge. The results are illustrated in Fig.~\ref{Fig_SINDYPI} (see Appendix for the data and method used). Further, these methods require complete knowledge of the underlying state space variables and data with higher sampling rates, which is rarely available in practice. 

\begin{figure}	
 \centering
\includegraphics[width=9cm]{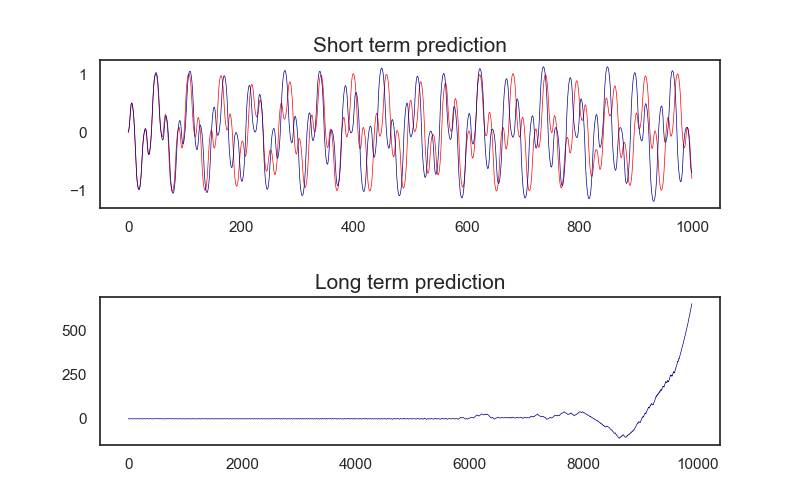}
\caption{A Double Pendulum Observation: In the top panel, a coordinate of the original time-series (in red) and the SINDY-PI predicted time-series (in blue) for 1000 time-steps are shown. All coordinates are used for training. In the bottom panel, the predicted values over 10000 time-steps is shown. The blow-up  clearly indicates the predicted time-series veering off from the set of interest.}
\label{Fig_SINDYPI}
\end{figure}

\begin{figure} 
 \centering
\includegraphics[width=9cm]{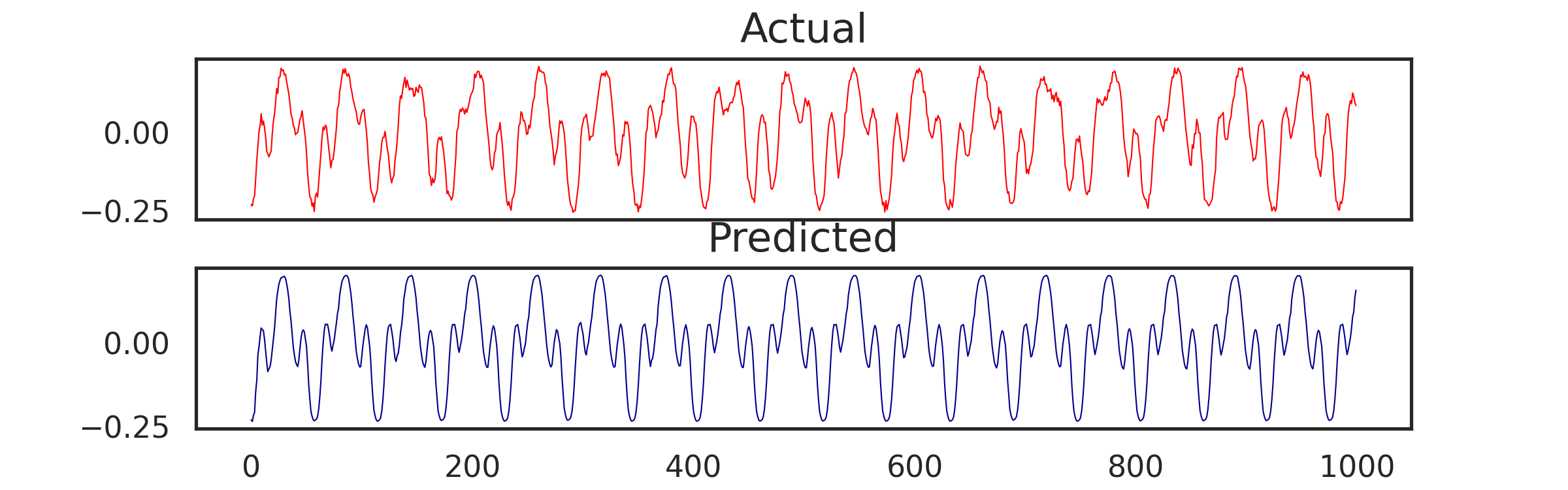}
      \caption{A Double Pendulum Observation: The map in $F_S$ in Takens Delay Embedding  (Top to Bottom) is learnt by by a deep neural network. First 1000 values of the actual time-series (in the top panel) and the predicted values (in the bottom panel) are shown.}
    \label{fig:Takens_DP_traj}
    \end{figure}

\begin{figure} 
 \centering
\includegraphics[width=9cm]{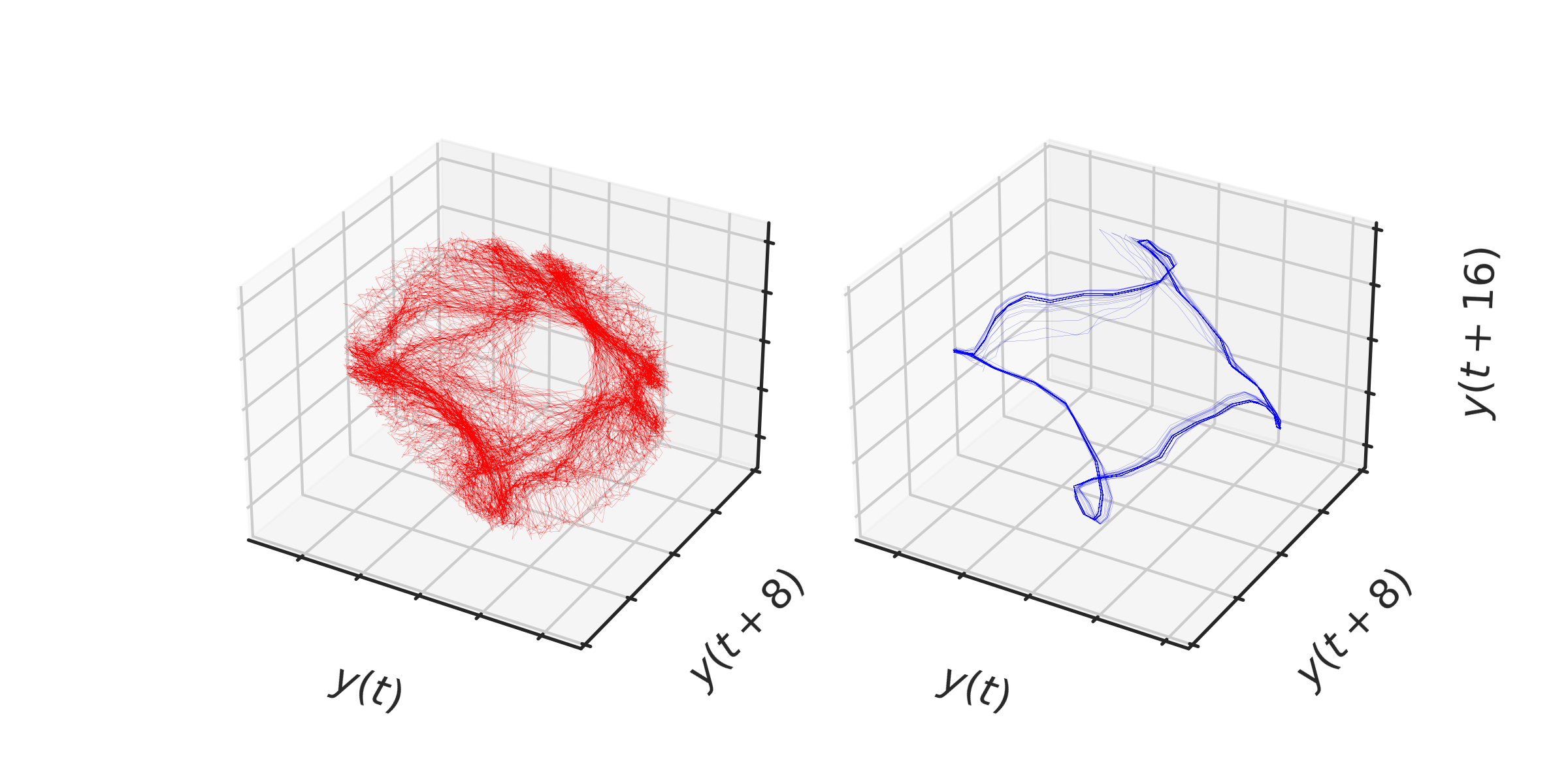}
      \caption{time-lag plot where three values $(y(t), y(t+8), y(t+16))$ are plotted for the two time-series in Fig.\ref{fig:Takens_DP_traj}}    \label{fig:Takens_DP_traj2}
    \end{figure}

Another aspect of Takens-delay embedding is that         global approximation strategies that discover a single data-fitting function $F_S$ frequently perform well, but often only if it has a \textit{low-functional complexity}, by which we mean a subjective measure of ``wiggliness" in the graph of a function.  The subtleness of the graph of the function can be missed even with deep learning techniques, and the resultant dynamics is very sensitive to it.    Takens along with his co-authors in \cite{bakker2000learning} allude to this issue when neural networks used for learning $F_S$  often fail (see  \cite{principe1992prediction}).  We demonstrate this through a numerical experiment, with a deep neural network. The actual values of the observations considered in Fig.~\ref{Fig_SINDYPI} and the time-series generated by employing $F_S$  are shown in  Fig.~\ref{fig:Takens_DP_traj}. Although the actual and reconstructed values look near-by, a closer look reveals that the predicted ones is close to a periodic orbit and misses several subtle features that are in the actual observation. The
 difference between the actual and forecasted dynamics is clearly evident when a three-dimensional recurrent or time-lag plot where three values $(y(t), y(t+8), y(t+16))$ are plotted in the three-dimensional coordinate space for the time-series shown in Fig.~\ref{fig:Takens_DP_traj2}.

Reservoir computing relies on processing the data using it as an exogenous input to a dynamical system of very large dimensions to achieve simpler training. 
We present the core mathematical ideas to the reader in Sec.~\ref{Sec_RC}.              
Existing machine learning methods in the reservoir computing framework \cite{jaeger2004harnessing, lu2018attractor} are more noise-resistant, have simpler training methods, and often choose to learn a ``linear" map that has low functional complexity -- the theory in \cite{grigoryeva2018echo} provides a mathematical basis for this approximation. However, the existence of a learnable map that can precisely predict the dynamics for a general recurrent neural network is unknown, and in practice results in experimenting with the reservoir. This learnability problem is solved only for a generic linear driven system in \cite{grigoryeva2023learning}, but the learnt system can suffer the same stability problem, where the forecasted trajectory can veer off while as in learning $F_S$ above. Methods such as \cite{gauthier2021next} that claim to avoid the need for a reservoir are special cases of Volterra filters \cite{boyd1985fading} completely fail on forecasting illustrated method (for sanity, the code is included in the repository).

A lack of learnability results in poor long-term consistency \cite{herteux2020breaking,Manju_IEEE, faranda2019boosting}.
  In summary, the current reservoir computing framework does not establish learnability for a given network while there is an inherent advantage of stability while considering nonlinear networks. We deal with the learnability problem by establishing the required topological conjugacy between the underlying system and the dynamics of a recurrent neural network.

\section{Realizing a SI-invertibility Condition with Reservoir Computers}
\label{Sec_RC}

In our theoretical framework, we consider a nonlinear driven system defined on a compact state space, and show how an adaptation of the customary with reservoir computing framework, one can present learnability, and also prevent the trajectories veering off.  In particular, we show the existence of an invertible function $f$ acting on the space $W$ considered above which has the property that $f(W)$ is contained in the global  attractor of a parametric dynamical system so that the trajectories do not wander far away from a desired set even when errors are made while learning.

When full coordinates of the underlying dynamical system are not available,
we first exploit Takens delay embedding theorem to generate an input sequence $\{u_i\}$ that represents the full coordinates of the state of a dynamical system, and then use two successive states of the reservoir $(x_{n-1}, x_{n})$ obtained from 
$(\ldots, u_{n-2}, u_{n-1})$ to forecast $u_n$, while in a customary reservoir computing framework only $x_n$ is used to forecast $u_n$.

The theory behind our approach is that we can create a one-to-one mapping between $(\ldots, u_{n-2}, u_{n-1})$ and the tuple $(x_{n-1}, x_{n})$ for all $n$. To make the learning feasible, we show that the mapping is continuous. To appeal to the dynamical systems community, we show that there is a dynamical system in the embedded space that is equivalent to the underlying system. Hence, in the theory presented our aim is 
to establish the entire commutativity ladder between $(W,S)$ and its equivalent system defined by $(Y_T, G_T)$ as in Fig.~\ref{Fig_COMM}~A. In the process, we also formally introduce the unfamiliar reader to the reservoir computing approach.

We consider the recursive dynamics of a continuous map $g:  U \times X 
\to X$. Here,  we assume  $U$ and $X$ to be nonempty compact metric spaces. We call a function $g$, a driven system with the input space $U$ and its ``state space" $X$ quietly understood.   A bi-infinite sequence $\{x_n\}$ is called a solution of $g$ for an input $\{u_n\}$ if   
$x_{n+1} = g(u_n,x_n)$ holds for all $n\in \mathbb{Z}$. Driven systems have at least one solution for a given input (e.g., \cite{manjunath2013echo}). A particularly interesting case:
if for each input $\bar{u} = \{u_n\}$ there exists exactly one solution, we say that $g$ has the unique solution property (USP) i.e., there is a well-defined solution map $\Psi : \bar{u} \mapsto \{x_n\}$. In the field of reservoir computing,  $\Psi(\bar{u})$ is meant to act as a proxy of the input sequence $\bar{u}$, but  
  we warn the reader that $\Psi$ is not necessarily injective when $g$ has the USP, for example, $g(u,x) := \frac{ux}{2}$regardless of any input sequence in $U$, $x_n\equiv 0$ for $n \in \mathbb{Z}$ is the only solution of $g$. However, $\Psi$  can be made injective with the SI-invertibility condition defined below.  The USP is equivalent to the echo state property \cite{manjunath2013echo} in the vast reservoir computing literature inspired by \cite{jaeger2004harnessing,maass2002real}, and is an asymptotic state contraction property that has been well-studied \cite{grigoryeva2018echo,manjunath2020stability}.     A remarkable feature of a driven system with the USP is that the ``future states"  of a solution need not be solved analytically, but  can be computed iteratively without the need for the entire-left infinite history of the input thanks to the ``uniform attraction property" (UAP) \cite{Manju_Nonlinearity}.  
Precisely, given any input $\bar{u}$ and its solution $\Psi(\bar{u})$, for every $\epsilon>0$, by initializing the driven system with an arbitrary initial value $y_m \in X$, then the ``state-trajectory" $y_{m+1}, y_{m+2},\ldots$ satisfying $y_{j+1}= g(u_j,y_j)$ for $j \geq m$ is such that the distance between  $x_{m+i}$ and $y_{m+i}$ is less than $\epsilon$ for all $i\ge k$ 
where (importantly) $k$ can be chosen independently of both $y_m$ and the input $\bar{u}$ (see \cite[Theorem 1]{Manju_Nonlinearity}).

We say $g$ satisfies the state-input (SI) invertibility condition, i.e,  $g(\cdot,x):U \to X$ is invertible for all $x$ (or in other words if $x_{n+1}=g(u_n,x_n)$ then $(x_{n+1},x_{n})$ uniquely determines $u_n$).  

{\bf Example of SI-invertibility.} We consider a recurrent neural network (RNN), a driven system with both $U, X \subset \mathbb{R}^N$ that is  defined through \begin{equation} \label{eq_RNNs}
	x_{n+1} = (1-a)x_n + a\overline{\tanh}(Au_n + \alpha Bx_n),
\end{equation}
where $A$ and $B$ are real invertible matrices of dimension $N \times N$ and $\overline{\tanh}(*)$ is the $\tanh$ performed component-wise on $*$, and $0 < \alpha$, $0<a\le 1$ are reals; a sufficient condition for $g$ to have the USP is that the norm of $\alpha B$ is less than $1$ (see \cite{Manju_ESP}). 
It suffices that when we allow $A$  and $B$ to be invertible ($A$ can be allowed to have a left-inverse),  $u_n$ can be recovered given $(x_{n}, x_{n+1})$ using 
\begin{equation} \label{eqn_SI_RNN}
u_n := A^{-1}\Big(\overline{\tanh}^{^{-1}}\frac{1}{a}\Big(x_{n+1}-(1-a)x_n\Big)  - \alpha B x_n \Big).
\end{equation}

\section{Embedding an Inverse-limit system of a dynamical system}

A dynamical system in this work is a tuple $(U,T)$ where $U$ is a compact metric space and $T: U \to U$ is a surjective function (that is not necessarily continuous).  
When $T$ is not surjective to start with we can restrict the non-transient dynamics to an invariant set $A$
(i.e., $T(A):=\cup_{u\in A} T(u) = A$) and $T$ restricted to $A$ is surjective. Due to the surjectivity, 
We call a bi-infinite sequence $\bar{u} = \{u_n\}_{n\in \mathbb{Z}}$ that obeys the update equation, $u_{n+1}=T(u_n)$ where $n \in \mathbb{Z}$ as an orbit/trajectory of $T$, and since $T$ is surjective there is at least one orbit ``through" each point in $U$.

Finding an equivalent dynamical system to $(U,T)$ means finding another dynamical system $(U',T')$ so that there exists a homeomorphism $\phi:U \to U'$ with the property that  $\phi \circ T=T'\circ \phi$ or equivalently the following diagram commutes:
\begin{equation}  \label{comm_h_repeat}
%    \[ 
    \psset{arrows=->, arrowinset=0.25, linewidth=0.6pt, nodesep=3pt, labelsep=2pt, rowsep=0.7cm, colsep = 1.1cm, shortput =tablr}
 \everypsbox{\scriptstyle}
 \begin{psmatrix}
U & U\\%
 U' & U'.
 %%%
 \ncline{1,1}{1,2}^{T} \ncline{1,1}{2,1} <{\phi}
 \ncline{1,2}{2,2} > {\phi}
 \ncline{2,1}{2,2}^{T'}
 \end{psmatrix}
% \]
\end{equation}
Such a map $\phi$ is called a  conjugacy and we say that $(U',T')$ is conjugate to $(U,T)$ or simply $T'$ is conjugate to $T$. If we relax the condition on $\phi$ where, instead of having a homeomorphism, we only require $\phi$ to be continuous, then we call $\phi$ a semi-conjugacy, and say that $T'$ is semi-conjugate to $T$.  When  $T'$ is semi-conjugate to $T$ then it is customary to say that $T$ is a (dynamical) extension of $T'$ and captures the dynamics of $T$, while $T'$ cannot have any additional complexity (measured through topological entropy) than present in $T$ (e.g., \cite{de2013elements}). Our aim here is to build a model that is either conjugate or a (dynamical) extension of the underlying system so that the underlying system can be forecasted.

Driven systems have an analog of 
semi-conjugacy too \cite{Manju_Nonlinearity}. For learning the underlying input dynamics using a driven system, its solution should not have additional ``complexity" (see \cite{Manju_Nonlinearity}) than that is present in the input. Hence there needs to be a function called the universal semi-conjugacy $h$ \cite{Manju_Nonlinearity} mapping a left-infinite sequence in $U$ to $X$ so that 
\begin{equation} \label{eqn_semi-conjugacy}
h \circ \sigma_v  = g(v,\cdot) \circ h,
\end{equation}
where $\sigma_v((\ldots, u_{-2},u_{-1}))=(\ldots, u_{-2},u_{-1},v)$.  The function $h$ here is also called a functional by other authors (e.g.\cite{boyd1985fading, grigoryeva2018echo}).

When the inputs to $g$ are restricted to be orbits of $(U,T)$ the states of the system at any time instant $n$ can be considered to depend on the left-infinite history of the input $(\ldots,u_{n-2},u_{n-1}$), i.e, on a left-infinite orbit of $T$.  The space of left-infinite orbits 
$$
 \widehat{U}_T: = \{(\ldots,u_{-2},u_{-1}) : Tu_{n} = u_{n+1} \},
$$ 
is called the inverse-limit space of $(U,T)$ (e.g., \cite{ingram2011inverse}); in the literature, the elements of $\widehat{U}_T$ are customarily written as right-infinite sequences. The map $T$ also induces a self-map $\widehat{T}$ on 
$\widehat{U}_T$ defined by  $\widehat{T}: (\ldots,u_{-2},u_{-1}) \mapsto 
(\ldots,u_{-2},u_{-1},T(u_{-1}))$.  It is straightforward to observe that $\varphi_T : U \to \widehat{U}_T$, $\varphi_T: u \mapsto (\ldots, u_{-2}, u_{-1}, u)$ (where $T(u_{-k}) = u_{-k+1}$ and $u_0=u$) is a conjugacy between $T$ and $\widehat{T}$, i.e.,
\begin{equation} \label{eqn_newconj}
	\widehat{T} \circ \varphi_T = \varphi_T \circ  T. 
\end{equation}

\begin{Theorem}{{\bf (Causal Embedding Theorem.)}}
 \label{Thm_CET}
	Let $g$ be a driven system that is SI-invertible and has the USP. Let $h$ denote the universal semi-conjugacy and $H_2(\overleftarrow{u}) := (h(r\overleftarrow{u}),h(
	\overleftarrow{u}))$,  where $r$ is the right-shift $r: (\cdots, u_{-2},u_{-1}) \mapsto 
  (\cdots, u_{-3},u_{-2})$.  Let $(\widehat{U}_T, \widehat{T})$  be the inverse-limit system of any dynamical system $(U,T)$. 
 The function $H_2$ restricted to $\widehat{U}_T$ is a topological semi-conjugacy between the inverse-limit system 
  $(\widehat{U}_T, \widehat{T})$  and the induced dynamical system  
  $(Y_T,G_T)$, i.e., the following diagram commutes
\begin{equation} \label{comm_H_CET}
\psset{arrows=->, arrowinset=0.25, linewidth=0.6pt, nodesep=3pt, labelsep=2pt, rowsep=0.7cm, colsep = 1.1cm, shortput =tablr}
 \everypsbox{\scriptstyle}
 \begin{psmatrix}
 \widehat{U}_T & \widehat{U}_T\\%
 Y_T &  Y_T.
 %%%
 \ncline{1,1}{1,2}^{\widehat{T}} \ncline{1,1}{2,1} <{H_2}
 \ncline{1,2}{2,2} > {H_2}
 \ncline{2,1}{2,2}^{G_T}
 \end{psmatrix}
 \end{equation}
or in other words, $(Y_T, G_T)$ is a factor of  $(\widehat{U}_T, \widehat{T})$. Further, if $T:U \to U$ is a homeomorphism, then $H_2$ embeds $\widehat{U}_T$ in $X_U \times X_U$, and hence
$(Y_T, G_T)$ is conjugate to $(\widehat{U}_T, \widehat{T})$.
\end{Theorem}

{\bf Proof. } Consider the inverse-limit space $\widehat{U}_T \subset \overleftarrow{U}$, and a bi-infinite orbit $\bar{u}$ of $T$. The left-infinite segment $\cev{u} = (\ldots,u_{-2},u_{-1})$ then belongs to $\widehat{U}_T$. If we restrict the map $\sigma_v$ to $\widehat{U}_T$, and select $v=Tu_{-1}$, then we have 
$\sigma_v(\ldots,u_{-2},u_{-1})=\sigma_{Tu_{-1}}(\ldots,u_{-2},u_{-1})=\widehat{T}(\ldots,u_{-2},u_{-1})$.

Let $(x_{n-1},x_n) \in Y_T$. Without loss of generality, let $n=0$, i.e., $(x_{-1},x_0) \in Y_T$. By definition of $Y_T$, there exists an orbit $\{u_k\}$ of $T$  for which $\{x_n\} = \Psi(\{u_k\})$ 
so that $(x_{-1},x_0) \in Y_T$ (recall $\Psi$ is the solution map). Also by \eqref{eqn_semi-conjugacy}, and since $T$ is surjective, $h(\cev{u}) = x_0$ and $h(r\cev{u}) = x_{-1}$ where $\cev{u} = (\ldots,u_{-2},u_{-1})$ is the left-infinite segment of the orbit $\{u_k\}$. 
Hence $H_2: \widehat{U}_T \to Y_T$ is surjective. 
 Since $h$ and $r$ are continuous, it follows that $H_2$ is also continuous. 
Also, since $h\circ \widehat{T}(\cev{u}) = x_1$, we have $H_2 \circ \widehat{T}(\cev{u}) =  (x_{0},x_{1}) = G_T(x_{-1},x_{0})$. Using all of these in \eqref{eqn_semi-conjugacy},  the diagram in \eqref{comm_H_CET} commutes, and $(Y_T,G_T)$ is a factor of $(\widehat{U}_T, \widehat{T})$. 

It remains to be shown whenever $T$ is a homeomorphism  then $H_2$ does not map two distinct points in $\widehat{U}_T$ to the same point. Suppose there exists two orbits $\{u_k\}$ and $\{v_k\}$ so that $\cev{u} \not= \cev{v}$ and
$H_2(\cev{u}) = H_2(\cev{v})$. This means that the solutions $\{x_k\} = \Psi(\{u_k\})$ and $\{y_k\}= \Psi(\{v_k\})$ are such that $(x_{-1},x_0) = 
(y_{-1},y_0)$. Since $g$ is SI-invertible, this implies $u_{-1}=v_{-1}$. Since $T$ is a homeomorphism, $u_{-1}$ and $v_{-1}$ have the same left-infinite history which implies $\cev{u} = \cev{v}$. Hence when $T$ is a homeomorphism, $H_2$ embeds $\widehat{U}_T$ in $X_U \times X_U$. By definition of $Y_T$, $H_2( \widehat{U}_T) = Y_T$ and hence $H_2: \widehat{U}_T \to Y_T$ is a homeomorphism. Thus $(Y_T,G_T)$ is topologically conjugate to $(\widehat{U}_T, \widehat{T})$. $\blacksquare$

We say a driven system $g$ \emph{causally embeds} a dynamical system $(U,T)$ if it satisfies two properties: (i) a universal semi-conjugacy $h$ exists and (ii)  a function $H_2: \widehat{U}_T \to X \times X$, defined by
$H_2((\ldots,u_{-2},u_{-1})) := (h((\ldots,u_{-3},u_{-2})),h((\ldots,u_{-2},u_{-1})))$ embeds the inverse-limit space $\widehat{U}_T$ in $X \times X$, i.e., there exists a set $Y_T \subset X \times X$ so that $H_2: \widehat{U}_T \to Y_T$ is a homeomorphism. The remarkable feature is that the ability of a driven system to causally embed a dynamical system does not depend on the map $T$ on $U$ since 
 $H_2$ is purely determined by $g$, and so the domain of $H_2$  -- from a practical viewpoint, this means that if the self-map $T$ is changed to another $\mathfrak{T}$ on $U$, the driven system $g$ does not need to change and when $H_2$ is considered or restricted to be defined only on $\widehat{U}_{\mathfrak{T}}$, it still embeds $\widehat{U}_{\mathfrak{T}}$ in $X \times X$. Hence, $H_2$ acts as a universal $X \times X$-valued  observable for any dynamical system on $U$ just like the delay-coordinate map $\Phi_{2d,\theta}$ is a universal observable 
 that does not depend on the map $S$ on $W$. Our interest in this universal observable is due to the map $h$ (and thus $H_2$ here) that can also be a differentiable conjugacy  since $h$ inherits the smoothness properties from $g$ \cite{grigoryeva2023learning}, and if $g_\lambda$ is continuous in the parameter $\lambda$, so would the conjugacy. With the causal embedding, we now show how there is a well-defined dynamical system defined on a subspace of $X \times X$ and a readout defined on that subspace.

\begin{Theorem} \label{Thm_GT_Exists}
Let $g$ be a driven system that is SI-invertible, and $(U,T)$ be a dynamical system. Consider the subspace $Y_T$ of $X_U \times X_U$ to be the tuple arising from two successive points of any solution of $g$, i.e.,
$$ \begin{multlined} Y_T:= \{(x_{n-1},x_n): \{x_k\}_{k\in \mathbb{Z}} \mbox{ is a solution for } \\  \mbox{some orbit of } T \mbox{ and } n \in \mathbb{Z}\}.  \end{multlined}$$ Then we have a well-defined map $G_T: Y_T \to Y_T$ defined by $G_T : (x_{n-1},x_n) \mapsto (x_n,x_{n+1})$. Consequently, the mapping $\Gamma: (x_{n-1},x_n) \mapsto u_{n}$ is well-defined when $x_{n-1}$ and $x_n$ are successive points on a solution obtained for an input 
$\{u_k\}$ that is an orbit of $T$.
\end{Theorem}

{\bf Proof.} Since $x_n=g(u_{n-1},x_{n-1})$ and $g(\cdot,x):U \to X$ is invertible, we have 
\begin{equation}
u_{n-1} = g_{*,x_{n-1}}^{-1}(x_n),
\end{equation}
where  $g_{*,x}^{-1}$ is the inverse of the map $g(\cdot,x):U \to X$. Using this in $x_{n+1}=g(u_n,x_n)=g(Tu_{n-1},x_n)$, we have $x_{n+1}=g(Tg_{*,x_{n-1}}^{-1}(x_n),x_n)$. Therefore, there exists a function $\theta_T: (x_{n-1},x_n) \mapsto x_{n+1}$. As a consequence we can define a map $G_T : (x_{n-1},x_n) \mapsto 
(x_n,x_{n+1})$ by 
\begin{equation} \label{eqn_GT}
	G_T: (x_{n-1},x_n) \mapsto 
(x_n,\theta_T(x_{n-1},x_n)). 
\end{equation}
Lastly, since $u_{n-1} = g_{*,x_{n-1}}^{-1}(x_n)$ and $u_{n} = Tu_{n-1}$, there exists a mapping $(x_{n-1},x_{n}) \mapsto u_n$.
$\blacksquare$

\section{Precise Forecasting due to the Conjugacies}

\begin{figure}	
%\flushleft
 \centering
\includegraphics[width=9cm]{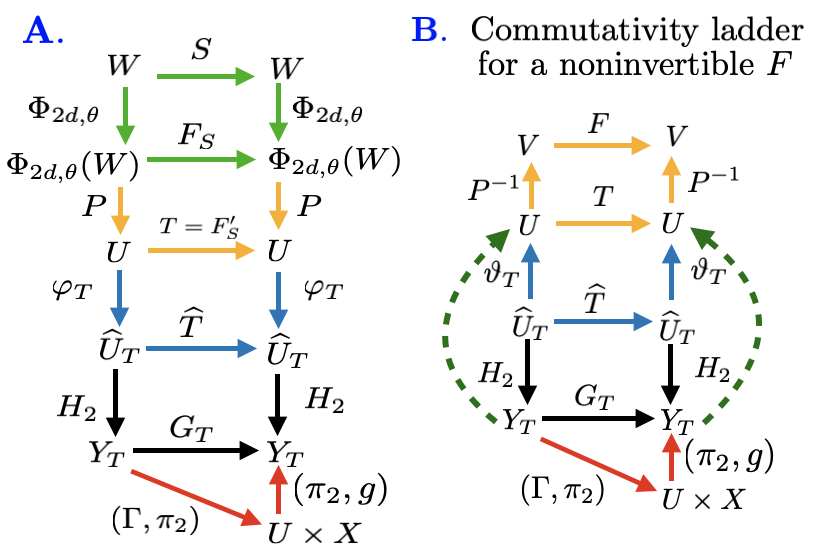}
\caption{\textcolor{blue}{A.} Commutativity ladder: inner workings of  a conjugate system when the Takens theorem holds. \textcolor{blue}{B.}  Commutativity ladder with input being an orbit from a non-invertible system.}
\label{Fig_COMM}
\end{figure}

In the customary RC approach (e.g., \cite{jaeger2004harnessing,lu2018attractor}), an RNN of the type \eqref{eq_RNNs} is used where a matrix  $A$ has a dimension $N \times K$ where $K<N$ and the observations $\{\theta(w_n)\}$ (where $w_{n+1} = S(w_n)$ and $\theta: W \to [-1,1]$) is fed into the RNN sequentially. Here we use a RNN in \eqref{eq_RNNs} where $A$ and $B$ have the same dimension $N$. We now produce arguments to establish equivalences of the dynamical systems that we have defined so far. For that we consider sever intermediary systems.  The reader may note that when we have a dynamical system $(Z_1,f_1)$ and there is a homeomorphism $r: Z_1 \to Z_2$, then $(Z_2,f_2)$ is conjugate to  $(Z_1,f_1)$ since 
$f_2 = r \circ f_1 \circ r^{-1}$.  The following maps and systems are considered:

1. $\Phi_{2d,\theta}: W \to \Phi_{2d,\theta}(W) \subset \mathbb{R}^{2d+1}$ -- the delay coordinate map induces $F_S$ conjugate to the underlying map $S$; see  Fig.~\ref{Fig_COMM}~A (printed in green) when Takens holds.

2. $P: \mathbb{R}^{2d+1} \mapsto U \subset \mathbb{R}^N$, where $P: v_n \mapsto (v^1_n,v^2_n,\ldots,v^{2d+1},0,0,\ldots 0)$ is the projection map. Since the projection is a homeomorphism onto its image,  the system $(U,T) :=  (P(\Phi_{2d,\theta}(W)), F_S')$ is conjugate to $F_S$. For the commutativity, see  Fig.~\ref{Fig_COMM}~A (printed in green-yellow)

3. When the RNN $g$ can causally embed $(U,T)$, then we have a map $G_T$ on $Y_T$ by Theorem~\ref{Thm_GT_Exists}, where $Y_T$ is the image of $\widehat{U}_T$ under $H_2$. For the commutativity, see  Fig.~\ref{Fig_COMM}~A (printed in blue-black).

The other conjugacy between $T$ and $\widehat{T}$ in  Fig.~\ref{Fig_COMM}~A shown in yellow-blue is due to \eqref{eqn_newconj}. In the portion of Fig.~\ref{Fig_COMM}~A printed in black-red, $\pi_2$ is defined as $\pi_2 : (a,b) \mapsto b$ regardless on which product space the tuple $(a,b)$ lies.

We can learn $\Gamma$ in the ``driven mode" of $g$ where we make a sufficient collection of tuples $(x_{n-1},x_{n})$ by feeding a finite segment of an orbit of $T$ after the dynamics locks onto a solution due to UAP. Finally, we realize the iterates of $G_T$ by: 
$G_T(x_{n-1},x_{n}) = (x_n, g(\Gamma(x_{n-1},x_{n}),x_{n}))$ and call the process of reading out $u_n=\Gamma(x_{n-1},x_{n})$  as an ``autonomous run" of $g$.  We note that the embedding $f = H_2 \circ \varphi_T \circ P \circ \Phi_{2d,\theta}$ of $W$ in $X \times X$ qualifies as a generalized synchronization in \cite{grigoryeva2021learning} where existential conditions are known for generic linear RNNs.

{\bf When Takens delay embedding does not hold.} In such a case, we presume that the inputs to $g$ are directly the elements of an orbit of an unknown dynamical system $(V,F)$; here $V \subset \mathbb{R}^K$ and we do not assume the surjective map  $F$ to be necessarily continuous or invertible, and we set $P:V \to U$ to be the zero-padding map defined previously to obtain a system $(U,T)$ whose orbits are vectors of length $N$ (see Fig.~\ref{Fig_COMM}~B in yellow).  The maps $G_T$ and  $\Gamma$ of course exist whenever $g$ is SI-invertible.

There is an advantage of learning $\Gamma$ rather than learning $F$ when both $T$ is invertible or non-invertible as explained in the following section.

\section{Global dissipativity (GD).} 
Learning error can be detrimental for forecasting.
Think about the logistic family of maps that has been studied extensively in the dynamical systems world. The family is defined by $T(u) = p \,u(1-u)$ on the real line, where $p$ is a real parameter. A particular map from this family when $p=4$ has the property that the image of the interval $[0,1]$ under $T$ is itself, or in other words $[0,1]$  is invariant. This full-logistic map is a classical example of chaotic behaviour. Suppose a machine learning algorithm is naively used to learn  $T$ and the learnt map is $T(u) = (4 +\epsilon) \,u(1-u)$, where $\epsilon> 0$, then except for a small set all other forward orbits would tend to  $- \infty$. RC can offer forecasting robustness, in the sense that a situation never arises where the dynamics veers off to an undesired set.

To understand the robustness of forecasting in the autonomous mode, let us suppose that due to a small amount of computational error, the map  $\Gamma$ we learn is not exact but has an error on a nonempty subset of $Y_T$. Due to this error in learning, the image of the readout $Y_T$ could be a space that is different from the intended $U$. Since the input is fed back to the reservoir in the autonomous run, the states of the reservoir may fall out of the set $Y_T$.  We would like to analyse the impact of the error on the stability of  the forecasted time-series in the autonomous run. 
For that purpose, let us assume that a readout  is found from a parameter family of readouts $P$, and let $P$ be a compact metric space.

Denote the readout $\Gamma^+$ here to be a continuous function defined on $P \times X \times X$. Note that we assume that it is defined on the whole of $P \times X \times X$ to take into account that the successive states could expand beyond $Y_T$. If $p_0 \in P$, then $\Gamma^+$ is the extension of the exact $\Gamma$ that we intended to learn on $Y_T$, i.e., when $\Gamma^+(p_0,\cdot,\star)$ restricted to $Y_T$ is the same as $\Gamma$. Hence the dynamical system  $G_{T,p_{0}}: (x,x') \mapsto 
(x', g(\Gamma^+(p_0,x,x'), x')$ has $Y_T$ as an invariant set (recall that for a dynamical systems $F:Z \to Z$, a set $A$ is invariant if $F(A)=A$).

\begin{figure}	
%\flushleft
 \centering
\includegraphics[width=8.5cm]{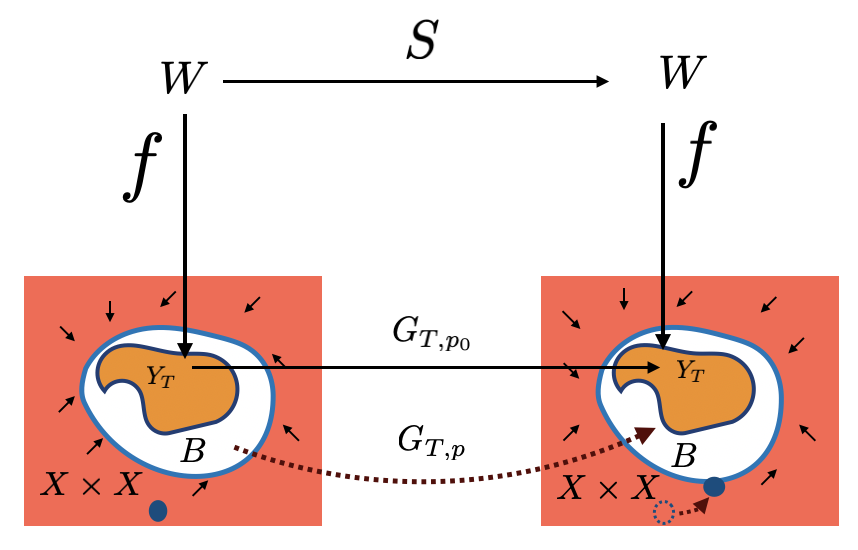}
\caption{Schematic of GD. Here $f = H_2 \circ \varphi_T \circ P \circ \Phi_{2d,\theta}$ of $W$ in $X \times X$ is the embedding. The ``cushion" around $Y_T$ is shown with the blue boundary.} 
\label{Fig_GDn}
\end{figure}

We say $g$ has the global dissipativity (GD) property if there exists a closed set $B\supset Y_T$ that contains all of the invariant sets of the parametric dynamical system defined by $G_{T,p}: (x,x') \mapsto 
(x', g(\Gamma^+(p,x,x'), x')$, where $p\in P$. The set $B$ then produces a ``cushion" around $Y_T$, since all the asymptotic dynamics of a dynamical system on a compact set is contained in an invariant set. Note that the successive states are always contained in the compact set $X \times X$, which is not necessarily invariant when $g(u,x)$ is contracting in the $x$-variable. What GD ensures is that we are considering a cushion $B$  
that does not leave out any asymptotic dynamics. 
When GD holds, the scenario is as depicted in Fig.~\ref{Fig_GDn}.

To show GD holds for $g$, we need to establish the existence of a closed set contained in $X \times X$ that is the union of the invariant sets of the parametric dynamical system defined by $G_{T,p}$ above, where $p\in P$. In general, when $F: Z \to Z$ is a continuous map on a metric space $Z$, the set $A := \cap_{n>0} F^n(Z)$ is the global attractor of $F$, and we here refer it as the attractor. Note that whenever $Z$ is non-empty and compact, the attractor always exists, and  $A$ is also non-empty and compact due to $Z$ being so.  Further, $A$ contains all invariant subsets of $F$ as well \cite{de2013elements}.

Now,  consider the autonomous system defined on $P \times X \times X$ defined by $
 F^+ : (p,x_{n-1},x_{n}) \mapsto (p, x_{n}, g(\Gamma^+(p,x_{n-1},x_{n}),x_{n}))$, where the dynamics in the $p$-variable is just the identity map. 

Let $A$ be the attractor of $F^+$. So $A$ contains all the invariant sets of $F^+$, and in particular contains the invariant set $\{p_0\} \times Y_T$. Since the dynamics in the $p$ coordinate is the identity map, the set $A$ can be written as $A = \bigcup_{p \in P} \{p\} \times A_p$, where $A_p$ is the attractor of $G_{T,p}$. Thus if we set  $B$ to be the union of all the attractors $A_p$ 
GD holds. Note that $A_{p_{0}}$  contains $Y_T$ and hence $B$ contains $Y_T$.

In summary, the cushion $B$ captures all the asymptotic dynamics  of the reservoir dynamics, and thereby the forecasted values do not veer off due to learning or training errors.

\section{Numerical Experiments} \label{Sec_NE}

\begin{center}
\begin{figure}[t]
 \centering
\includegraphics[width=8.5cm]{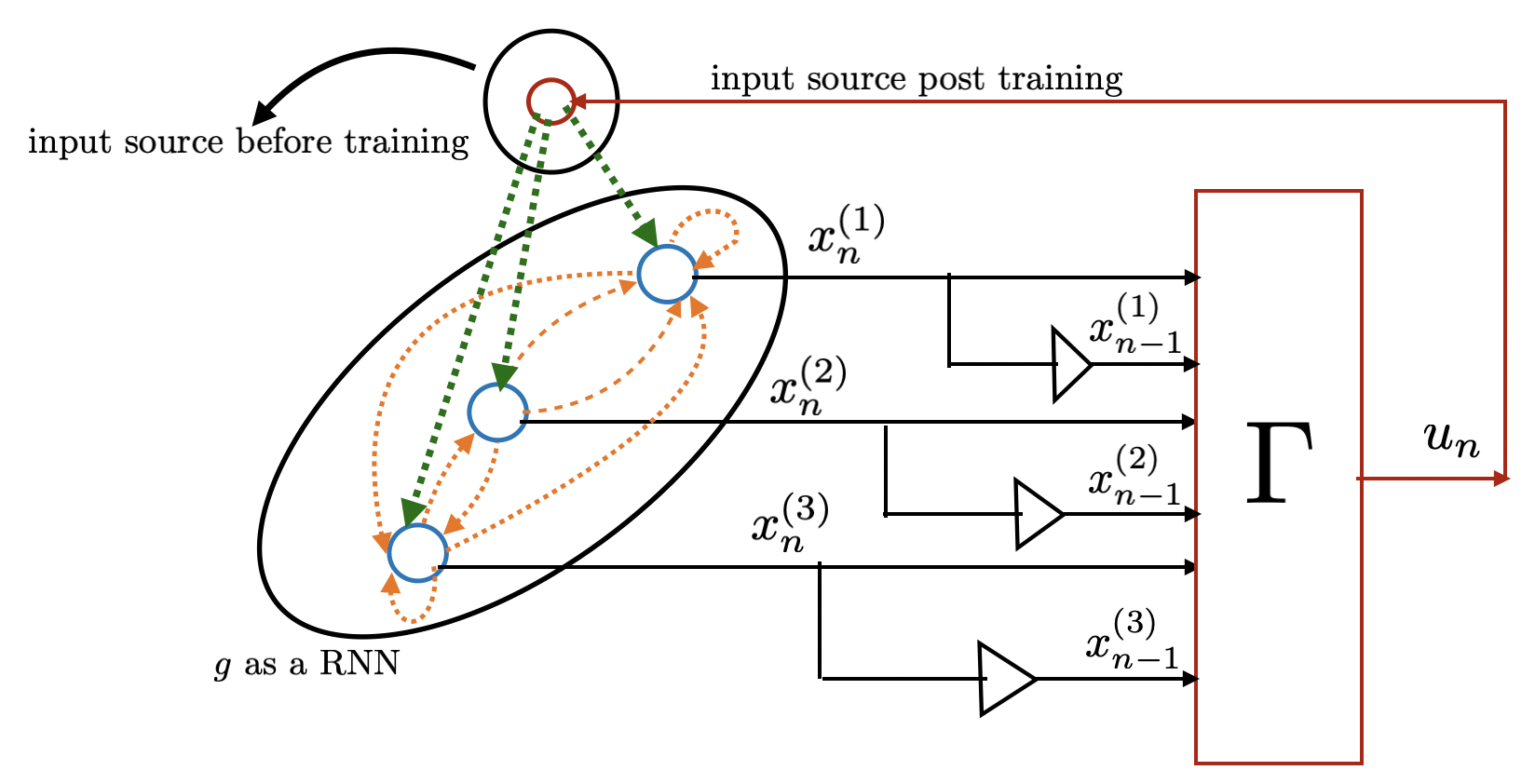}
\caption{Schematic of a RCN with  three nodes in the  RNN: (the print in red indicates the connections made to run in the autonomous mode)}
\label{Fig_RCN}
 \end{figure}
\end{center}

The schematic post-training with a recurrent neural network (RNN) is shown in Fig.~\ref{Fig_RCN}, and due to the conjugacies, we call the system as a recurrent-conjugate network. 

 We consider dynamical systems where the limiting behavior of a time-series $\{\theta(w_n)\}$ can be observed, a feature of systems where $S$ is surjective. For example, the double pendulum (DP) system has a highly intricately woven invariant set owing to the fact that its dynamics is isomorphic to the product of a Bernoulli shift and a rotation of the circle \cite{ornstein1989ergodic}. 
 The inherent effect of a rotation on the circle is that the auto-correlation of a trajectory could be aperiodic instead of vanishing rapidly with the lag, a feature of many chaotic systems considered in the literature.  Building a model from data would require capturing the effect of the rotation and the chaotic dynamics of the Bernoulli shift. We consider the system of differential equations \eqref{DP_Eqn1m} and \eqref{DP_Eqn2m} (see Appendix) and data obtained from it for our experiment.  We describe the complete details of learning and training of $\Gamma$ using a deep neural network in the Appendix.

 \begin{figure}	
 \centering
\includegraphics[width=8.5cm]{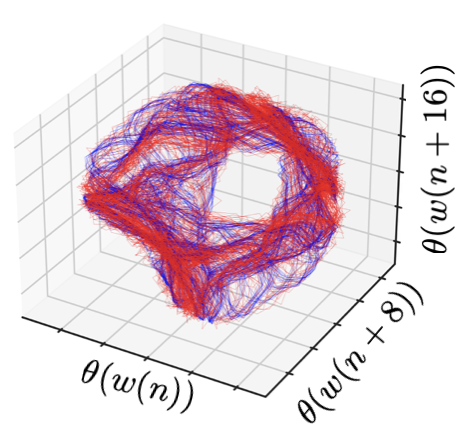}
\caption{Multi-lag recurrence plot of a forecasted time-series from an observation $\theta$ of a Double Pendulum System with a zero-mean normally distributed noise  with standard deviation  0.01; the points 
$(\theta(w(t)), \theta(w(t+8)), \theta(w(t+16)))$  of the actual (in red) and forecasted (in blue) values are plotted.} 
\label{Fig_Att_new}
\end{figure}

 \begin{figure}	
 \centering
\includegraphics[width=8.5cm]{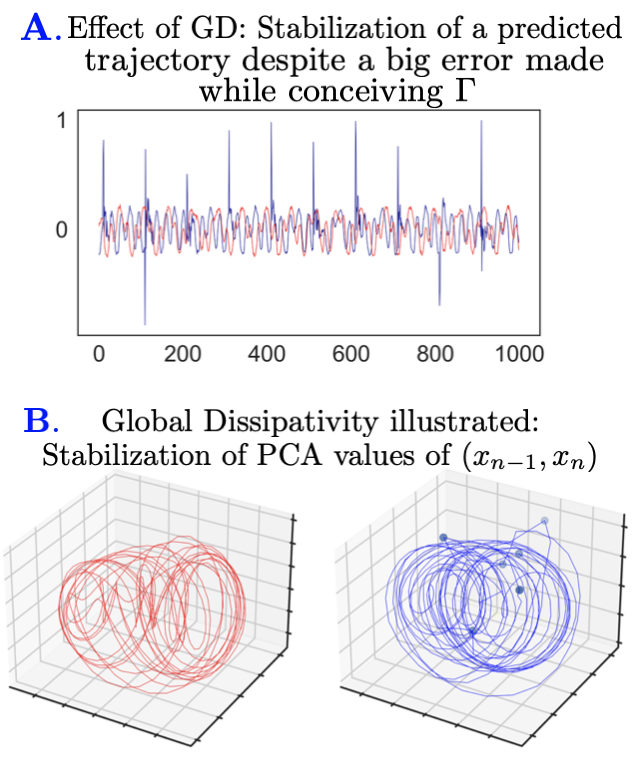}
\caption{\textcolor{blue}{A.} Effect of GD: The underlying tuples $(x_{n-1},x_{n})$  of the data in (E) is perturbed by a high-amplitude  Gaussian noise vector with standard deviation 2  once in every 100 time-steps in the autonomous mode of $g$ to simulate an error in $\Gamma$; despite this the readout in (blue) stabilizes towards the pattern of the actual readout conceived in the driven mode (red).
        \textcolor{blue}{B.}
      Illustration of GD: three largest principal component analysis (PCA) values of $(x_ {n-1},x_{n})$ corresponding to the time-series in the top panel is plotted for the driven (red) and autonomous (blue) mode; successive points are joined; thick dots show the high-amplitude noise instants. The blue trajectory drifts but is stabilized back. }
        \label{Fig_GD_new}
\end{figure}

  In Fig.~\ref{Fig_Att_new}, we present the  multi-lag recurrence plot $(\theta(w(t)), \theta(w(t+8)), \theta(w(t+16)))$ of both the actual and forecasted values. To illustrate the advantage of global dissipativity, a zero-mean normally distributed noise vector is added to $(x_{n-1},x_{n})$, (i.e., to an orbit of $G_T$) every 100 time-steps in the autonomous mode to simulate an excursion outside of $Y_T$. The  readouts when the system is driven (driven mode) in Fig.~\ref{Fig_GD_new} are compared with such a noise-introduced autonomous run. Further, due to the high dimension of $X$, we plot three leading principal component analysis (PCA) values of $(x_{n-1},x_{n})$ in the same figure to illustrate how a trajectory in the autonomous mode (shown in blue) is stabilized back. Both these pictures illustrate the remarkable effect of GD which enable RCNs to forecast systems like the double pendulum that has not been achieved before.   
  
  By contrasting the outcomes of other various popularly known methods, we further demonstrate the superiority of predicting using RCNs on other systems in \cite{Supp}.  A future study would identify conditions on the driven system for causal embedding to hold when an input $\{\theta(w_n)\}$ is fed sequentially rather than as delay coordinates, eliminating the requirement to know the length of the delay-coordinate vector. In \cite{Supp}, we present 
numerical findings by 
 feeding  $\{\theta(w_n)\}$ sequentially as well to address a real world problem.

\noindent {\bf Functional complexity of $G_T$.} 
Since we have used two successive states of the reservoir for training, and by learning $\Gamma$ we indirectly learn $G_T$, we 
we remark on how the functional complexity of $G_T$ varies as the reservoir size increases. More details are available in \cite{Supp}. As the dimension of $X$ i.e., $N$ gets larger than the effective dimension of the input delay-coordinate vector $K\ge 2d+1$ for a $g$ in \eqref{eq_RNNs}, empirical evidence suggests that the resultant map $G_T$ turns ``more linear" or has lower functional complexity as evidenced through the generalized Pearson correlation coefficients tending to 1 (see Table~\ref{Table_FC2}). We refer the reader to \cite[Table 1]{Supp} for more  details).

\begin{center}
\begin{table}
\centering
    \scalebox{0.85}{\begin{tabular}{|c|c| c c |} 
        \hline
        Input, from & Dimension of $X$ & $u_n$ vs $u_{n+1}$ 
        & $\begin{bmatrix}
            x_{n-1}\\
            x_n
        \end{bmatrix}$ vs $\begin{bmatrix}
            x_n\\
            x_{n+1}
        \end{bmatrix}$ \\
        \hline
        \hline
  Double Pendulum: & & & \\
         & 10 & 0.9710  & 0.9822\\       
          & 100 &  & 0.9803\\
         & 1000 & & 0.9854\\
\hline
       Lorenz map & & &\\
        & 10 &0.8568  & 0.9444\\
            & 100 &  & 0.9592\\
           & 1000 & & 0.9703\\
            
        \hline 
    H\'enon map: & & & \\
        & 10 &-0.3652  & 0.3524\\
        
        & 100 & & 0.7416\\
        & &  0.8529\\        \hline 
       Pomeau-Manneville  & & & \\
         & 10 &0.8743  & 0.9722\\        
            & 100 & & 0.9890\\
            & 1000 & &  0.9930\\
        \hline 
        Logistic map & & & \\
  & 10 &-0.0291  & 0.5615\\          
            & 100 &  & 0.7191\\
          & 1000 & & 0.7884\\ \hline
    \end{tabular}}
   \caption{Multidimensional Correlation Coefficients $\rho$ to indicate a general reduction in the functional complexity of the map $G_T$ that arises from a RNN. Each rows corresponds to a different dynamical system. The second column corresponds to the dimension of $X$ (here, number of artificial neurons) in the RNN used. 
    The numerical estimate of $\rho$ for the relevant vectors is plotted in the last two columns. The reader can compare the linear relationship $u_n$ vs $T(u_{n})$ and the linear relationship $\begin{bmatrix}
            x_{n-1}\\
            x_n
        \end{bmatrix}$ vs $G_T \left( \begin{bmatrix}
            x_{n-1}\\
            x_n
        \end{bmatrix}\right )$. The dimension of the RNN is also given, showcasing that increasing the dimension of the driven system typically indicates a map $G_T$ of lower functional complexity.} \label{Table_FC2}
    \end{table}
\end{center}

\section{Conclusions}

Mathematical models can attain many different forms ranging from very simple ones to the complex ones.
Very simple ones include  phenomenological or ``toy models, with some conceptual resemblance to the original system but no evident direct derivation from it. Less complex ones include interpretable models with fewer equations -- these can be derived either out of intuition or from approximating a vector field or maps. Many simple or simpler models can not only grossly under characterize the qualitative behaviour of the system but produce systems without an attractor or dissipative dynamics thereby producing dynamics that blows up. At the most complex level, models are equivalent to a full first-principles scientific description of all the intricacies involved in the entire problem. Such systems can have dozens or even hundreds of artificial neurons, and hence equations defining different aspects of the problem.  The advances in machine learning  methods together with the understanding of dynamical systems has made building models with numerous such equations possible to regenerate the dynamics of the underlying systems with increased precision and stability. 

In this work, we provided a Takens delay embedding's stability challenge by learning a driven system's single-delay dynamics. Employing a network in contrast to obtaining vector fields from data forfeits the interpretability of the learnt model,  but the inscrutability of our approach comes with a prize of learning systems that have failed with existing methods due to either a lack of learnability or stability.

\noindent {\bf Data and Code Availability. The data and code are available in \url{ https://github.com/MJSteynberg/ForecastingThroughCausalEmbedding}.

\noindent {\bf Acknowledgements. } The authors gratefully acknowledge Franco du Plessis for assistance with the SINDY-PI simulation. The authors also thank National Research Foundation, South Africa for financial support under Grant UID 115940.

\section{Appendix: Data Generation and Implementation Details}

We consider the system of ordinary differential equations (ODEs) that describe the angular acceleration of the masses of a double pendulum:   
{\footnotesize
\begin{equation}
	 \begin{matrix}
		\ddot{\phi} = \left[\frac{m_2 g}{\ell_1(m_1+m_2)}  \sin \psi \cos(\psi-\phi) + \frac{m_2}{m_1+m_2}\dot{\phi}^2\sin(\psi-\phi)\cos(\psi-\phi) \right. \\
		
		\qquad \left.-\frac{g}{\ell_1}\sin\phi+\frac{m_2 \ell_2}{\ell_1(m_1+m_2)}\dot{\psi}^2\sin(\psi-\phi) \right] \left\{1-\frac{m_2}{m_1+m_2}\cos^2(\psi-\phi)\right\}^{-1}
	\end{matrix}
		\label{DP_Eqn1m}
\end{equation}

\begin{equation}
	 \begin{matrix}
		\ddot{\psi} =\left[\frac{g}{\ell_2} \sin \phi \cos(\psi-\phi) - \frac{m_2}{m_1+m_2} \dot{\psi}^2\sin(\psi-\phi)\cos(\psi-\phi) \right. \\
		
		\qquad \left. -\frac{g}{\ell_2}\sin\psi - \frac{\ell_1}{\ell_2}\dot{\phi}^2\sin(\psi-\phi)\right] \left\{1-\frac{m_2}{m_1+m_2}\cos^2(\psi-\phi)\right\}^{-1},
	\end{matrix}
		\label{DP_Eqn2m}
\end{equation}}
with $m_i$ denoting the masses, $g$ the gravitational constant and $l_i$ the length of the two rods; we arbitrarily set $(m_1,m_2,\ell_1,\ell_2) = (0.75,2.5,1.15,1)$, with mass measured in kilograms and length in metres.  
The two equations are then converted to a first order system with four variables where the variables denote the angular velocity and angular position. We sample the flow of this Double Pendulum system by using numerical integration techniques (implemented using the {\it SciPy} library's ODE solver \cite{2020SciPy-NMeth}) using a step size of 0.05.  The data generated is used in three numerical experiments.

First, we refer to the results in Fig.~\ref{Fig_SINDYPI} where the data is used.  The SINDY-PI code in \cite{kaheman2020sindy} is adapted to the data generated by the samples, and the forecasting results are shown in Fig.~\ref{Fig_SINDYPI}, and the ODE generated is exhibited in \cite{Supp}. 

In the second experiment, the code used in \cite{gauthier2021next}. The system failed, and 
the code is included in the repository mentioned in the data availability statement.

For the RCN experiment in Section~\ref{Sec_NE}, the forecasting is done as in practical situation where an observation of the data is used along with noise. The observable (chosen arbitrarily) is the sum of the position of the $y$ coordinate of the first mass and the $x$-coordinate of the second mass. 
 While RCNs are used, forecasting commences from the very next time step where training is stopped. 
 
 {\bf Setting up the network. } We consider systems that exhibit recurrence where trajectories that are seen during training belong to their $\omega$-limit set.  We realize $g$ through a RNN of the type \eqref{eq_RNNs}} with the causal embedding property. Regardless of how the map $\Gamma$ is learnt and implemented, we call a driven system with a RNN that is SI-invertible and has the USP  a \emph{recurrent conjugate network} (RCN).  A schematic of a RCN is shown in Fig.~\ref{Fig_RCN}. Throughout, the matrices $A$ and $B$ in the RNN are randomly initialized, and in particular, we set the spectral radius of the matrix $\alpha B$ to be $\alpha$ by using a matrix $B$ with spectral radius  $1$. For $\alpha$ sufficiently small, which is usually in the interval $(0,1)$, RNNs have the USP \cite{yildiz2012re}; we make a passing remark that the upper permissible value of $\alpha$ for the USP to hold can also be verified empirically (using the parameter-stability plot in \cite{manjunath2020stability}).  We retain the same parameter values  $a=0.5$ and $\alpha=0.99$ for all the experiments. The random matrices are 
initialized with a seed so that the same network is used across simulations. In all simulations, 
the first 1000 datapoints were discarded to allow the network to `forget' its initial states or more precisely for the UAP to take effect so that the states are practical purposes are the solution values. More details are below:

\begin{itemize}
\item Step size for discretization = 0.05.
\item Observable $\theta(w_n)$: $y$-coordinate of the first mass + $x$-coordinate of the second mass (arbitrarily chosen) is rescaled so that 
$\theta(w_n)$ lies in $[-0.2,0.2]$ and has zero mean
\item input 
$$u_n:= \begin{bmatrix}
    \theta(w_n) \\
    \theta(w_{n-1})  \\
    \vdots \\
    \theta(w_{n-d+1})
\end{bmatrix} + 0.01z_n,$$  
where $z_n$ is zero-mean iid noise that is normally distributed with standard deviation 0.01.
\item $d=10$.
\item Washout = 1000 samples
	\item Length of Data used for training of $\Gamma = 5000$
	\item Length of Training data for learning $\Gamma$ = 4000.
	\item Length of Validation for learning $\Gamma$  = 1000.
	\item Length of Autonomous Run = 10000. 
\end{itemize}

\noindent {\bf RCN Implementation and Training.} In our RCN implementation, we learn $\Gamma: (x_n,x_{n-1}) \mapsto u_n$ by first obtaining the principal components of the states $x_n$ and then 
learning a feedforward network (NN) on these principal components with the target $u_n$. The principal components are used only for an efficient state representation possibly to reduce the errors while learning, is {\bf optional},  and is not for a lossy approximation since all principal components are used.  To be explicit, we denote the matrix with the first $N$ states of the network data as row vectors by $X_{1:N}$. If $X_{1:N}=U\Sigma P^T$ denotes the singular value decomposition of $X_{1:N}$, then the principal component matrix is given by $P$ and the principal components are given by 
 $Z_{1:N}=X_{1:N}P.$
These principal components are used to train a feedforward neural network. If we denote the row vectors of $Z_{1:N}$ by 
$z_i^T, i=1,2,\cdots,N$ and the neural network by $NN$, then we learn an approximation of the map
$NN: (z_{n-1},z_n) \mapsto  u_n.
$
We use the learnt approximation of $NN$ to approximate $\Gamma$ since 

{\footnotesize
$$
NN\left( \begin{bmatrix} 
P^Tx_{n-1} \\
P^Tx_n
\end{bmatrix}
\right) = NN \circ 
\begin{bmatrix}
P^T & 0 \\
0 & P^T 
\end{bmatrix}\begin{bmatrix}
x_{n-1}\\
x_n
\end{bmatrix} \\ = u_n = \Gamma(x_{n-1},x_n).
$$}

The feedforward neural network is implemented in {\it Python} using the {\it Keras} \cite{chollet2015keras} library built on {\it Tensorflow}. 

The network is constructed with 12 hidden layers with a layer dimension equal to 64. The activation function on the hidden layers is the {\it ReLU } function built into {\it Keras}, whereas the output layer has a {\it tanh} output. Training is accomplished using the {\it Adam optimizer}, optimizing the {\it Mean Squared Error} loss function. We train the network three times, using three different learning rates for the {\it Adam optimizer} equal to 0.001, 0.0001, and 0.00001. Each time we learn $\Gamma$, we use 150 training epochs and a batch size of 128.

\bibliographystyle{wmaainf}
\bibliography{pnas}

\end{document}